\documentclass[12pt]{amsart}
\usepackage{amssymb,latexsym,comment,url}

\newcommand{\PP}{{\mathbb P}}
\newcommand{\Q}{{\mathbb Q}}

\newcommand{\Z}{{\mathbb Z}}

\newenvironment{Proof}{\par\noindent{\sc Proof:}}%
                      {\hspace*{\fill}\nobreak$\Box$\par\medskip}
                       {\hspace*{\fill}\nobreak$\Box$\par\medskip}

\newtheorem{Proposition}{Proposition}[section]
\newtheorem{Theorem}[Proposition]{Theorem}
\newtheorem{Lemma}[Proposition]{Lemma}
\newtheorem{Corollary}[Proposition]{Corollary}
\newtheorem{Example}[Proposition]{Example}

\theoremstyle{definition}

\newtheorem{Remark}[Proposition]{Remark}

\renewcommand{\baselinestretch}{1.1}

\begin{document}

\title[On rational periodic points of $ x^d+c$]%
{On rational periodic points of $ x^d+c$}

\author[M. Sadek]%
{Mohammad~Sadek}
\address{American University in Cairo, Mathematics and Actuarial Science Department, AUC Avenue, New Cairo, Egypt}
\email{mmsadek@aucegypt.edu}
\let\thefootnote\relax\footnote{Mathematics Subject Classification: 37P05, 37P15 }
\begin{abstract}
   We consider the polynomials $\displaystyle f(x)=x^d+c$, where $d\ge 2$ and $c\in\mathbb Q$. It is conjectured that if $d=2$, then $f$ has no rational periodic point of exact period $N\ge 4$. In this note, fixing some integer $d\ge 2$, we show that the density of such polynomials with a rational periodic point of any period among all polynomials $f(x)=x^d+c$, $c\in\Q$, is zero. Furthermore, we establish the connection between polynomials $f$ with periodic points and two arithmetic sequences. This yields necessary conditions that must be satisfied by $c$ and $d$ in order for the polynomial $f$ to possess a rational periodic point of exact period $N$, and a lower bound on the number of primitive prime divisors in the critical orbit of $f$ when such a rational periodic point exists. The note also introduces new results on the irreducibility of iterates of $f$.
\end{abstract}
\maketitle

\section{Introduction}
An arithmetic dynamical system over a number field $K$ consists of a rational function $f: \PP^n(K)\to \PP^n(K)$ of degree at least $2$ with coefficients in $K$ where the $n^{th}$ iterate of $f$ is defined recursively by $f^{1}(x)=f(x)$ and $f^{m}(x)=f(f^{m-1}(x))$ when $m\ge 2$. A point $P\in\PP^n( K)$ is said to be a periodic (preperiodic) point for $f$ if the orbit $P,f(P), f^2(P),\cdots,f^n(P),\cdots$ of $P$ is periodic (eventually periodic). If $N$ is the smallest positive integer such that $f^N(P)=P$, then the periodic point $P$ is said to be of exact period $N$.

The following conjecture was proposed by Morton and Silverman. There exists a bound $B(D, n, d)$ such that if $K/\Q$ is a number field of degree $D$, and $f : \PP^n(K) \to \PP^n(K)$ is a morphism of degree $d\ge 2$ defined over $K$, then the number of $K$-rational preperiodic points of $f$ is bounded by $B(D,n,d)$, see \cite{Morton}. When $f$ is taken to be a quadratic polynomial over $\Q$, the following conjecture was suggested in \cite{Poonen}. If $N\ge 4$, then there is no quadratic polynomial $f(x)\in\Q[x]$ with a
rational point of exact period $N$. The conjecture has been proved when $N=4$, see \cite{Morton2}, and $N=5$, see \cite{Flynn}. A conditional proof for the case $N=6$ was given in \cite{Stoll2}.

We consider the polynomial $f(x)=x^d+c$ over a number field $K$. If $c=c_1/c_2$ where $c_1$ and $c_2$ are relatively prime in the ring of integers $O_K$ of $K$, we investigate the divisibility of the coefficients of the iterates $f^m(x)$, $m\ge 2$, by the prime divisors of $c_1$ and $c_2$. Using these divisibility criteria, we approach three questions concerning the arithmetic dynamical system of $f(x)=x^d+c$: (i) When is $f(x)$ stable over $K$? (ii) Fixing $d$, what is the density of such polynomials with periodic points? (iii) Given that $f(x)$ possesses a rational periodic point of period $n$, should this yield necessary conditions satisfied by $d$ and $c$?

The stability question in arithmetic dynamical systems concerns the irreducibility of the iterates of $f(x)$ over $K$. More precisely, a polynomial $f(x)$ is said to be stable over a field $K$ if $f^n(x)$ is irreducible over $K$ for every $n\ge 1$.
In \cite{Ahmadi}, the authors showed that most monic quadratic polynomials in $\Z[x]$ are stable over $\Q$. One may find sufficient conditions for an irreducible monic quadratic polynomial in $\Z[x]$ to be stable over $\Q$ in \cite{Jones1}. It was shown that $f(x)=x^2+c\in \Z[x]$ is stable over $\Q$ if $f(x)$ is irreducible itself, see \cite{Stoll}.
 Further, the polynomial $f(x)=x^d+c\in\Z[x]$, $d\ge 2$, is known to be stable over $\Q$ if $f(x)$ is irreducible, see \cite{DANIELSON}.

Unlike the situation over $O_K$, $f(x)=x^d+c\in K[x]$ can be irreducible over $K$ whereas $f^{n}(x)$ is reducible over $K$ for some $n>1$.
In this note, if $c=c_1/c_2$ where $c_1$ and $c_2$ are relatively prime in $O_K$, we show that the existence of a prime divisor $p$ of $c_1$ such that $\gcd(\nu_p(c_1),d)=1$, where $\nu_p$ is the valuation of $K$ at the prime $p$, implies the stability of $f(x)$. For instance, if $d$ is prime and $c_1$ is not a $d^{th}$-power modulo units in $O_K$, then $f(x)$ is stable.

 Assuming that $u_1/u_2$ is a periodic point of $f(x)$ of exact period $n$, where $u_1$ and $u_2$ are relatively prime in $O_K$, we give several results on the divisibility of the coefficients of the iterate $f^n(x)$ by prime divisors of $u_1$ and $u_2$. This enables us to show that if $f(x)$ has a $K$-rational periodic point, then $c_2$ must be a $d$-th power modulo units in $O_K$. More precisely, $c_2=u_2^d$ modulo units. Fixing $d$, a hight argument, then, yields that the density of such polynomials with periodic points among all polynomials $f(x)=x^d+c$ is zero. In particular, almost all polynomials $f(x)=x^d+c$ satisfy the conjecture of Morton and Silverman.

 We establish the connection between a periodic point $u_1/u_2$ of $f(x)=x^d+c\in\Q[x]$ of period $n$ and the sequence $u_1^{m}-u_2^{m}$, $m=1,2,\cdots$. In fact, we show that $c_1$ divides $u_1^{d^n-1}-u_2^{d^n-1}$, yet none of the prime divisors of $c_1$ divide $u_1-u_2$. This provides us with necessary conditions on $c_1$ in order for $f(x)$ to have such a periodic point. For instance, one knows that if $p$ is a prime divisor of $c_1$ such that $\gcd(p-1,d^n-1)=1$, then $f(x)$ has no periodic points of period $n$.

 Finally, we display the relation between rational periodic points of the polynomials $f(x)=x^d+c\in\Q[x]$ and another sequence, namely the sequence of the iterates, $f^n(0)$, evaluated at $0$. One may consult \cite{Krieger} for several results on the existence of primitive prime divisors of such sequences. In this note, we show that the existence of a periodic point of $f(x)$ of exact period $n$ implies a lower bound on the number of primitive prime divisors of $f^n(0)$.

\section{Valuations of the coefficients of the iterates of $f$}
In this section, we assume that $K$ is an arbitrary field unless otherwise stated.
\begin{Lemma}
\label{lem:freecoefficient}
Let $f(x)=x^d+c$, $d\ge 2$, $c\in K$. One has $f^{n}(0)=c+c^d g_n(c)$ where $g_n\in \Z[x]$ is a polynomial of degree $d^{n-1}-d$, $n\ge2$.
\end{Lemma}
\begin{Proof}
Since $\displaystyle f^2(0)=c+c^d$, the statement is true when $n=2$ by taking $g_2(x)=1$. Now, an induction argument will yield the statement. Assume that $f^{n}(0)=c+c^dg_n(c)$ where $g_n(x)\in \Z[x]$ is of degree $d^{n-1}-d$. One has that $f^{n+1}(0)=c+(f^n(0))^d$. One observes that
\begin{eqnarray*}
f^{n+1}(0)&=&c+(c+c^dg_n(c))^d=c+c^d\left(1+c^{d-1}g_n(c)\right)^d.
\end{eqnarray*}
We set $g_{n+1}(x)=\left(1+x^{d-1}g_n(x)\right)^d$.
 The polynomial $g_{n+1}(x)\in \Z[x]$. Moreover, since $g_n$ has degree $d^{n-1}-d$ by assumption, one gets that the degree of $g_{n+1}$ is $d(d^{n-1}-d+d-1)=d^{n}-d$.
\end{Proof}
The following lemma gives an explicit description of the coefficients of $f^n(x)$.
\begin{Proposition}
\label{prop:valuations}
Let $f(x)=x^d+c$, $d\ge 2$, $c\in K$. Assume that $f^n(x)=f_0+f_1x^d+f_2x^{2d}+\ldots+f_{d^{n-1}}x^{d^n}$. The following statements are correct.
\begin{itemize}
\item[a)] $f_{d^{n-1}}=1$.
\item[b)] $f_i\in c \Z[c]$ for every $0\le i< d^{n-1}$.
\item[c)] $\deg f_i= d^{n-1}-i$ for $0\le i\le d^{n-1}$.
\end{itemize}
\end{Proposition}
\begin{Proof}
That $f_0\in c\Z[c]$ and $\deg f_0$ in $\Z[c]$ is $d^{n-1}$ is implied by Lemma \ref{lem:freecoefficient}.

We now follow an induction argument. For the polynomial $f^2(x)$, one has
\begin{eqnarray*}
f^2(x)=(x^d+c)^d+c&=&x^{d^2}+\sum_{i=0}^{d-1}  {d\choose i} c^{d-i} x^{id} +c\\
&=& x^{d^2}+c \sum_{i=1}^{d-1}  {d\choose i} c^{d-1-i} x^{id} +c+c^d.
\end{eqnarray*}
Since $\displaystyle f_i={d\choose i} c^{d-i}\in c\Z[c]$, $1\le i< d-1$, is of degree $<d$, the statement is correct for $f^2(x)$.

Assume the statement holds for $f^n(x)$. One obtains the following equalities
\begin{eqnarray*}
f^{n+1}(x)&=&(f^n(x))^d+c=\left[f_0+f_1x^d+f_2x^{2d}+\ldots+f_{d^{n-1}-1}x^{d^{n}-d}+x^{d^n}\right]^d+c\\
&=&\left[c\left(f'_0+f'_1x^d+f'_2x^{2d}+\ldots+f'_{d^{n-1}-1}x^{d^{n}-d}\right)+x^{d^n}\right]^d+c
\end{eqnarray*}
where $f'_i =f_i/c\in \Z[c]$ and $\deg f'_i<d^{n-1}-1$ by assumption. Setting $f'(x)= f'_0+f'_1x^d+f'_2x^{2d}+\ldots+f'_{d^{n-1}-1}x^{d^{n-1}}$, one obtains
\begin{eqnarray*}
f^{n+1}(x)&=&x^{d^{n+1}}+\sum_{i=1}^d{d\choose i}c^if'(x)^ix^{d^n(d-i)}+c.
\end{eqnarray*}
It is obvious that each coefficient of $f^{n+1}(x)-x^{d^{n+1}}$ is in $c\Z[c]$.

For part c), one sees that \[f^{n+1}(x)=(f^n(x))^d+c=\left(f_0+f_1x^d+f_2x^{2d}+\ldots+f_{d^{n-1}-1}x^{d^{n}-d}+x^{d^n}\right)^d+c.\] We are looking for the degree of the coefficient of $x^{ld}$ in the latter expansion where $0\le l\le d^{n}$. Using an induction argument, we assume that $\deg f_i=d^{n-1}-i$ in $\Z[c]$. In view of the multinomial expansion, the latter expansion is given by
\[f^{n+1}(x)=\sum_{k_0+k_1+\ldots+k_{d^n}=d}{d\choose{k_0,\ldots, k_{d^n}}}\prod_{t=0}^{d^n}(f_tx^{td})^{k_t}+c.\]
Using the induction assumption, the degree of the coefficient of $x^{ld}$ in $f^{n+1}(x)$ is obtained as follows
\begin{eqnarray*}
\sum_{t=0}^{d^n}k_t(d^{n-1}-t)&=&d^{n-1}\sum_{t=0}^{d^n} k_t-\sum_{t=0}^{d^n}tk_t\\
\end{eqnarray*}
where $\displaystyle \sum_{t=0}^{d^n}k_t=d$ and $\displaystyle \sum_{t=0}^{d^n}tk_td=ld$.
\end{Proof}
The following corollary is a straight forward result of the proposition above.
\begin{Corollary}
\label{cor1}
Let $K$ be a discrete valuation field with ring of integers $O_K$.
Let $f(x)=x^d+c$, $d\ge 2$, where $c=c_1/c_2$ is such that $c_1$ and $c_2$ are relatively prime in $O_K$. Assume that $f^n(x)=f_0+f_1x^d+f_2x^{2d}+\ldots+f_{d^{n-1}}x^{d^n}$. Then
\[c_2^{d^{n-1}}f^n(x)=F_0(c_1,c_2)+F_1(c_1,c_2)x^d+F_2(c_1,c_2)x^{2d}+\ldots+F_{d^{n-1}-1}(c_1,c_2)x^{d^{n}-d}+F_{d^{n-1}}(c_1,c_2)x^{d^n}\]
where $F_i(c_1,c_2)=c_2^{d^{n-1}}f_i\in \Z[c_1,c_2]$ is a homogeneous polynomial of degree $d^{n-1}$. Moreover, $F_i(c_1,c_2)\in c_1c_2^i\Z[c_1,c_2]$ if $i\ne d^{n-1}$; and $F_{d^{n-1}}(c_1,c_2)=c_2^{d^{n-1}}$ .
\end{Corollary}
\begin{Proof}
 Since $f_i \in c\Z[c]$, $i\ne d^{n-1}$, and $\deg f_i =d^{n-1}-i$ for $0\le i\le d^{n-1}$, see Proposition \ref{prop:valuations}, we may clear the denominators of the coefficients $f_i$'s by multiplying throughout by $c_2^{d^{n-1}}$, hence the result is obtained.
\end{Proof}
\section{The stability of $f(x)=x^d+c$}
 Let $K$ be a field with valuation $\nu$ whose value group is $\Z$. Let $F[x]\in K[x]$ be the polynomial $F_0+F_1x+\ldots+F_kx^k$ where $F_0\ne 0$ and $F_k\ne 0$.

The Newton polygon of $F$ over $K$ is constructed as follows. We consider the following points in the real plane: $A_i=(i,\nu(F_i))$ for $i=0,\ldots, k$. If $F_i=0$ for some $i$, then we omit the corresponding point $A_i$. The {\em Newton polygon} of $F$ over $K$ is defined to be the lower convex hull of these points. More precisely, we consider the broken line $P_0P_1\ldots P_l$ where $P_0=A_0$, $P_1=A_{i_1}$ where $i_1$ is the largest integer such that there are no points $A_i$ below the line segment $P_0P_1$. Similarly, $P_2$ is $A_{i_2}$ where $i_2$ is the largest integer such that there are no point $A_i$ below the line segment $P_1P_2$. In a similar fashion, we may define $P_i$, $i=2,\ldots,l$, where $P_l=A_k$. If some line segments of the broken line $P_0P_1\ldots P_l$ pass through points in the plane with integer coordinates, then such
points in the plane will be also considered as vertices of the broken line. Therefore, we may add $s\ge 0$ more points to the vertices $P_0P_1\ldots P_l$. The Newton polygon of $F$ over $K$ is the polygon $Q_0Q_1\ldots Q_{l+s}$ obtained after relabelling all these points from left to the right, where $Q_0=P_0$ and $Q_{l+s}=P_l$.

The following theorem generalizes Eisenstein's criterion of irreducibility, see for example \cite[Theorem 9.1.13]{Cheze}.
\begin{Theorem}[Eisenstein-Dumas Criterion]
\label{thm:Eisenstein}
Let $K$ be a field with valuation $\nu$ whose value group is $\Z$. Let $F(x)=F_0+F_1x+\ldots+F_kx^k\in K[x]$ with $F_0F_k\ne 0$. If the Newton polygon of $F$ over $K$ consists of the only line segment from $(0, m)$ to $(k, 0)$
and if $\gcd(k, m)=1$, then $F$ is irreducible over $K$.
\end{Theorem}

 We recall that $x^d+c$ is irreducible over a field $K$ if and only if for every prime $p$ dividing $d$, $-c$ is not a $p^{th}$-power in $K$; and if $4\mid d$ then $c$ is not $4$ times a $4^{th}$-power in $K$, see \cite[Theorem 8.1.6]{Karpilovsky}.

\begin{Theorem}
\label{thm1}
Let $K$ be a number field with ring of integers $O_K$. Let $f(x)=x^d+c$, $d\ge 2$, be such that $c=c_1/c_2$ is such that $c_1$ and $c_2$ are relatively prime in $O_K$. Assume that there is a prime $p$ in $O_K$ such that $\gcd(\nu_p(c_1),d)=1$ where $\nu_p$ is the valuation of $K$ at the prime $p$. Then $f(x)$ is stable over $K$.
\end{Theorem}
\begin{Proof}
Let $K_p$ be the completion of $K$ with respect to the prime $p$ and $\nu_p$ be the corresponding valuation. In view of Corollary \ref{cor1}, one has $\displaystyle f^n(x)=\frac{H_n(x)}{c_2^{d^{n-1}}}$ where \[H_n(x)=F_0(c_1,c_2)+F_1(c_1,c_2)x^d+F_2(c_1,c_2)x^{2d}+\ldots+F_{d^{n-1}-1}(c_1,c_2)x^{d^{n}-d}+F_{d^{n-1}}(c_1,c_2)x^{d^n}\] and $F_i(c_1,c_2)=c_2^{d^{n-1}}f_i$. Now we consider the Newton polygon of the polynomial $H_n(x)\in \Z[c_1,c_2][x]$ over $K_p$. According to Lemma \ref{lem:freecoefficient}, one has $\nu_p(F_0(c_1,c_2))=\nu_p(c_1)$. Proposition \ref{prop:valuations} indicates that $\nu_p(F_i(c_1,c_2))\ge \nu_p(c_1)$ if $1\le i< d^{n}$ and $\nu_p(F_{d^n}(c_1,c_2))=\nu_p\left(c_2^{d^{n-1}}\right)=0$ where the latter equality follows from the fact that $c_1$ and $c_2$ are relatively prime. Therefore, the Newton polygon of $H_n(x)$ consists of one line segment joining the two points $(0,\nu_p(c_1))$ and $(d^n,0)$. Since $\gcd(\nu_p(c_1),d^n)=1$ by assumption, Theorem \ref{thm:Eisenstein} yields that $H_n(x)$ is irreducible over $K_p$, hence over $K$. This implies that $f(x)$ is stable.
\end{Proof}
\begin{Corollary}
Let $K$ be a number field and $f(x)=x^d+c$, $d\ge 2$, where $c=c_1/c_2$ is such that $c_1$ and $c_2$ are relatively prime in the ring of integers $O_K$ of $K$. Assume that $c_1$ is not of the form $u v^p$ for any prime divisor $p$ of $d$, where $v\in O_K$ and $u$ is a unit of $O_K$. Then $f(x)$ is stable over $K$.\\
In particular, if $f(x)=x^d+c$ where $d$ is prime, then $f(x)$ is stable over $K$ if $c_1$ is not a $d^{th}$-power modulo units in $O_K$.
\end{Corollary}

In what follows, we see some examples of polynomials $f(x)$ violating the relative primality condition $\gcd(\nu_p(c_1),d)=1$ in Theorem \ref{thm1}. We remark that these polynomials are not stable.
\begin{Example}
  If one considers the polynomial $f(x)=x^d-c^d$, $c\in K$, over a field $K$, then $f(x)$ is not stable as $f^1(x)=f(x)$ is reducible.
The polynomial $f(x)=x^2-4/3$ is irreducible over $\Q$ since $4/3$ is not a square in $\Q$, yet $\displaystyle f^2(x)=\left(x^2-2x+\frac{2}{3}\right)\left(x^2+2x+\frac{2}{3}\right)$.
\end{Example}
\section{Periodic points}
From now on $K$ is a number field with ring of integers $O_K$. We will write $O_K^{\times}$ for the group of units in $O_K$. If $p$ is a prime in $O_K$, then $\nu_p$ is the valuation of $K$ at $p$.

We consider $f(x)=x^d+c$ where $c=c_1/c_2$ such that $c_1\in O_K$ and $c_2\in O_K/ O_K^{\times}$ are relatively prime in $O_K$. Given $u\in K$, the orbit of $u$ under $f$ is the set $\displaystyle O_f(u)=\left\{u,f(u),f^2(u),\ldots\right\}$. By a periodic point $u$ of exact period $n$, we mean that $f^n(u)=u$ and that $n$ is the smallest such positive integer. In particular, the polynomial $f^n(x)-x$ has a zero at $u$ and $O_f(u)$ is a finite set with exactly $n$ elements. Moreover, any point in the orbit $O_f(u)$ is a periodic point with period $n$. In particular, $f^n(x)-x$ has at least $n$ linear factors.

In accordance with Corollary \ref{cor1}, one recalls that
\begin{eqnarray*}
f^n(x)=\frac{F_0(c_1,c_2)+F_1(c_1,c_2)x^d+F_2(c_1,c_2)x^{2d}+\ldots+F_{d^{n-1}-1}(c_1,c_2)x^{d^{n}-d}+F_{d^{n-1}}(c_1,c_2)x^{d^n}}{c_2^{d^{n-1}}}.
\end{eqnarray*}
Finding the zeros of $f^{n}(x)-x$ is equivalent to finding the zeros of the following polynomial
\begin{eqnarray*}
\label{eq2}
G^n(x)=F_0(c_1,c_2)-c_2^{d^{n-1}}x+F_1(c_1,c_2)x^d+F_2(c_1,c_2)x^{2d}+\ldots+F_{d^{n-1}-1}(c_1,c_2)x^{d^{n}-d}+F_{d^{n-1}}(c_1,c_2)x^{d^n}.
\end{eqnarray*}
Given that $u_1/u_2$ is a periodic point of period $n$ of $f(x)$, where $u_1$ and $u_2$ are relatively prime in $O_K$ and $u_2\in O_K/ O_K^{\times}$, one multiplies throughout times $u_2^{d^n}$ to get
\begin{eqnarray}\label{eq1}
F_0u_2^{d^n}-c_2^{d^{n-1}}u_1u_2^{d^n-1}+F_1u_1^du_2^{d^n-d}+F_2u_1^{2d}u_2^{d^n-2d}+\ldots+F_{d^{n-1}-1}u_1^{d^{n}-d}u_2^d+F_{d^{n-1}}u_1^{d^n}=0\nonumber\\
\end{eqnarray}
where $F_i:=F_i(c_1,c_2)$.
\subsection{The denominators $c_2$ and $u_2$ of $c$ and $u$}
\begin{Proposition}
\label{prop:denominators}
Let $f(x)=x^d+c_1/c_2$ such that $c_1\in O_K$ and $c_2\in O_K/ O_K^{\times}$ are relatively prime in $O_K$. Let $u_1/u_2$ be a periodic point of $f(x)$ with period $n$ where $u_1,u_2\in O_K$ are relatively prime. The following properties hold.
\begin{itemize}
\item[a)] $u_2^d\mid F_{d^{n-1}}=c_2^{d^{n-1}}$.
\item[b)] $c_2$ and $F_0$ are relatively prime in $O_K$.
\item[c)] $c_2\mid u_2^{d^n}$.
\item[d)] $c_2$ and $u_2$ have exactly the same prime divisors.
\end{itemize}
\end{Proposition}
\begin{Proof}
(a) follows directly from eq (\ref{eq1}) and the fact that $u_1$ and $u_2$ are relatively prime in $O_K$.

For (b), Lemma \ref{lem:freecoefficient} yields that \begin{eqnarray*}F_0&=&c_1c_2^{d^{n-1}-1}+c_2^{d^{n-1}-d}c_1^dg_n(c_1/c_2),\qquad g_n(x)=\sum_{i=0}^{d^{n-1}-d}g_{n,i}x^i, g_{n,i}\in\Z\\
&=&  c_1c_2^{d^{n-1}-1}+c_2^{d^{n-1}-d}c_1^d\sum_{i=0}^{d^{n-1}-d}g_{n,i}(c_1/c_2)^i\\
 &=& c_1c_2^{d^{n-1}-1}+\sum_{i=0}^{d^{n-1}-d}g_{n,i}c_1^{d+i}c_2^{d^{n-1}-d-i}
\in c_1\Z[c_1,c_2].\end{eqnarray*}
Every term in the latter expansion of $F_0$ is divisible by $c_2$ except for the term whose coefficient is $g_{n,d^{n-1}-d}=1$. Since $c_1$ and $c_2$ are relatively prime, it follows that $c_2$ and $F_0$ are relatively prime in $O_K$.

For (c), since $F_i\in c_2^i\Z[c_1,c_2]$ except when $i=0$, see Corollary \ref{cor1}, this yields that $c_2\mid F_0 u_2^{d^n}$, see eq (\ref{eq1}). Since by (c), one knows that $c_2$ and $F_0$ are relatively prime, it follows that $c_2\mid u_2^{d^n}$. Part (d) follows from (a) and (c).
\end{Proof}
\begin{Corollary}
\label{cor:Z}
Let $c\in O_K$. If $f(x)=x^d+c$, $d\ge 2$, has a periodic point $u$, then $u\in O_K$.
\end{Corollary}
\begin{Proof}
This follows from Proposition \ref{prop:denominators} (d).
\end{Proof}
\begin{Theorem}
\label{Thm:denominators}
Let $f(x)=x^d+c_1/c_2$, $d\ge 2$, such that $c_1\in O_K$ and $c_2\in O_K/ O_K^{\times}$ are relatively prime in $O_K$. Let $u_1/u_2$ be a periodic point of $f(x)$ where $u_1,u_2\in O_K$ are relatively prime. One has $c_2= u_2^d$.
\end{Theorem}
\begin{Proof}
We assume that $u_1/u_2$ is of period $n$.
Let $p$ be a prime divisor of $u_2$. Proposition \ref{prop:denominators} d) implies that $p$ divides $c_2$. Considering eq (\ref{eq1}), one sets $\alpha:=\nu_p(c_2^{d^{n-1}}u_1u_2^{d^n-1})=d^{n-1}\nu_p(c_2)+(d^{n}-1)\nu_p(u_2)$. We also set
\begin{eqnarray*}
\alpha_l:&=& \nu_p(F_lu_1^{ld}u_2^{d^n-ld})=\nu_p(F_l)+(d^n-ld)\nu_p(u_2),\;0< l< d^{n-1}\\
&\ge& l\nu_p(c_2)+(d^n-ld)\nu_p(u_2)\\&=&d^n\nu_p(u_2)+l(\nu_p(c_2)-d\nu_p(u_2)),
\end{eqnarray*}
see Corollary \ref{cor1}.
Furthermore, we define \[\alpha_{d^{n-1}}:= \nu_p(F_{d^{n-1}})=\nu_p(c_2^{d^{n-1}})=d^{n-1}\nu_p(c_2),\qquad\alpha_0:=\nu_p(F_0u_2^{d^n})=d^n\nu_p(u_2),\]
see Corollary \ref{cor1} and Proposition \ref{prop:denominators} b), respectively.

If $\nu_p(c_2)<d\nu_p(u_2)$, then
\[\min_{0\le l< d^{n-1}}\alpha_l> d^{n-1}\nu_p(c_2)=\alpha_{d^{n-1}}.\]  In this case, either $\alpha_{d^{n-1}}=\alpha_r$ for some $r \ne d^{n-1}$, which is impossible, or $\alpha_{d^{n-1}}=\alpha$ which is again impossible as $\nu_p(u_2)>0$.

If $\nu_p(c_2)>d\nu_p(u_2)$, then \[\min_{0< l\le d^{n-1}}\alpha_l> d^{n}\nu_p(u_2)=\alpha_0.\] In the latter case, since $\alpha_0\ne \alpha_r$ for any $r\ne 0$, one must have $\alpha_0=\alpha$. It follows that $\nu_p(u_2)=d^{n-1}\nu_p(c_2)$ which contradicts our assumption that $\nu_p(c_2)>d\nu_p(u_2)$.

One concludes that it must be the case that $\nu_p(c_2)=d\nu_p(u_2)$ for any common prime divisor of $c_2$ and $u_2$. Therefore, assuming that $u_2\in O_K/O_K^{\times}$, one obtains that $c_2=u_2^d$.
\end{Proof}

\begin{Remark}
\label{rem1}
  If $u_1/u_2$ is a periodic point of $f(x)=x^d+c_1/c_2$ where $c_i$ and $u_i$ are as in Theorem \ref{Thm:denominators}, then $c_2=u_2^d$. In other words, a periodic point of $f(x)$ of any period will have the same denominator. In particular, if $f^j(u_1/u_2)=v_{1,j}/v_{2,j}$, $j=1,2,\ldots$, are the elements in the orbit $O_f(u_1/u_2)$ of $u_1/u_2$, where $v_{1,j}$ and $v_{2,j}$ are relatively prime in $O_K$, then one may assume that $v_{2,j}=u_2$ for every $j$. In fact, since $c_2=u_2^d$, one has $f(u_1/u_2)=(u_1^d+c_1)/u_2^d$. Therefore, $u_2^{d-1}\mid (u_1^d+c_1)$.
\end{Remark}
The following is a direct consequence of Theorem \ref{Thm:denominators}.
\begin{Corollary}
\label{cor:families with no periodic points} If $f(x)=x^d+c_1/c_2$, $d\ge 2$, where $c_1$ and $c_2$ are relatively prime and $c_2$ is not a $d^{th}$-power in $O_K$, then $f$ has no periodic points of any period. In particular, there are infinitely many polynomials $f(x)=x^d+c$ that have no periodic points of any period.
\end{Corollary}
\begin{Corollary}
\label{cor:conjecture}
 Let $u_1/u_2$ be a periodic point of exact period $n$ of $f(x)=x^d+c_1/c_2$, where $c_i$ and $u_i$ are as above. If $g(x)=x/u_2^{d-1}$ and $h(x)=x^d+c_1$, then $u_1$ is a periodic point of the polynomial $g\circ h\in K[x]$ of exact period $n$.
\end{Corollary}
\begin{Proof}
  Recall that since $c_2=u_2^d$, see Theorem \ref{Thm:denominators}, one has $f(u_1/u_2)=(u_1^d+c_1)/u_2^d$. As $f(u_1/u_2)$ is an element in $O_f(u_1/u_2)$, it follows that $u_2^{d-1}$ divides $u_1^d+c_1$, see Remark \ref{rem1}. In other words, $f(u_1/u_2)=(g\circ h(u_1))/u_2$, where $g\circ h(u_1),u_2\in O_K$ are relatively prime. Now the statement follows by a simple induction argument to show that $f^j(u_1/u_2)=(g\circ h)^j(u_1)/u_2$. Now the statement of the corollary holds because $f^n(u_1/u_2)=u_1/u_2$.
\end{Proof}
Corollary \ref{cor:families with no periodic points} can be strengthened in the following manner over $\Q$. We recall that for $c=a/b\in \Q$ where $\gcd(a,b)=1$, one may define the {\em height} of $c$ to be $\displaystyle h(c)=\max\{|a|,|b|\}.$ Fixing $d\ge 2$, we define the following two subsets in $\Q$
\begin{eqnarray*}
 S(N)&=&\left\{\frac{\alpha}{\beta}:\alpha\in\Z,\,\beta\in\Z^{+},\,\gcd(\alpha,\beta)=1,\,h\left(\frac{\alpha}{\beta}\right)\le N\right\},\\
 S_d(N)&=&\left\{\frac{\alpha}{\beta}:\alpha\in\Z,\,\beta\in\Z^{+},\,h\left(\frac{\alpha}{\beta}\right)\le N,\, \beta \textrm{ is a $d$-th power}\right\}.
 \end{eqnarray*}
We will show that $\displaystyle \lim_{N\to\infty}\frac{|S_d(N)|}{|S(N)|}=0$. This implies the following consequence. Fixing $d\ge 2$, if $f(x)=x^d+c_1/c_2\in\Q[x]$, where $c_1\in\Z$ and $c_2\in\Z^+$ are relatively prime in $\Z$, has a periodic point, then $c_2$ is a $d$-th power. In other words, if we consider the set of such polynomials with periodic points such that the height of $c_1/c_2$ is less than $N$, then according to Theorem \ref{Thm:denominators}, the set of those $c_1/c_2$ is contained in $S_d(N)$. This means that the density of polynomials $x^d+c$ which have periodic points among all polynomials of the form $x^d+c$, $c\in\Q$, is zero. This can be restated as follows: Fixing $d\ge 2$, almost all polynomials $x^d+c$, $c\in\Q$, have no periodic points.
\begin{Proposition}
\label{prop:asymptotic}
For an integer $d\ge 2$, one has the following asymptotic formula \[\displaystyle \frac{|S_d(N)|}{|S(N)|}\sim\frac{\pi^2}{6 N^{(d-1)/d}}\quad\textrm{ as }N\to \infty.\]
\end{Proposition}
\begin{Proof}
It is clear that $|S_d(N)|$ is asymptotically $2N^{(d+1)/d}$. A standard analytic number theory exercise shows that $$\displaystyle\sum_{0<\alpha,\beta\le N,\,\gcd(\alpha,\beta)=1} 1$$ is asymptotically $6N^2/\pi^2$. It follows that $|S_d(N)|/|S(N)|$ is asymptotically $\displaystyle \frac{2\pi^2}{12 N^{(d-1)/d}}$.
\end{Proof}
Fixing $d\ge 2$, we set
\begin{eqnarray*}
P(N)&=& \left\{ c\in \Q: h(c)\le N\right\},\\
  P_d(N)&=& \left\{ c\in \Q: x^d+c \textrm{ has a periodic point},\; h(c) \le N\right\}.
\end{eqnarray*}
According to Theorem \ref{Thm:denominators}, one has $|P_d(N)|/|P(N)|<|S_d(N)|/|S(N)|$. Now, the following result holds as a direct consequence of Proposition \ref{prop:asymptotic}.
\begin{Theorem}
One has the following limit $\displaystyle \lim_{N\to\infty}\frac{P_d(N)}{P(N)}=0$.
\end{Theorem}
The above limit holds if one replaces $\Q$ with a number field. The proof is similar but the hight function has to be changed appropriately.
\subsection{The numerators $c_1$ and $u_1$ of $c$ and $u$}
We now deduce some divisibility conditions on the numerators of $c$ and $u$.
Recall that
\begin{eqnarray*}
G^n(x)=F_0-c_2^{d^{n-1}}x+F_1x^d+F_2x^{2d}+\ldots+F_{d^{n-1}-1}x^{d^{n}-d}+F_{d^{n-1}}x^{d^n},
\end{eqnarray*} and eq (\ref{eq1}) is given by
\begin{eqnarray*}
F_0u_2^{d^n}-c_2^{d^{n-1}}u_1u_2^{d^n-1}+F_1u_1^du_2^{d^n-d}+F_2u_1^{2d}u_2^{d^n-2d}+\ldots+F_{d^{n-1}-1}u_1^{d^{n}-d}u_2^d+F_{d^{n-1}}u_1^{d^n}=0.\nonumber\\
\end{eqnarray*}
In the following lemma, we list some of the divisibility criteria satisfied by the numerator $u_1$ of a periodic point $u_1/u_2$ of $f(x)=x^d+c_1/c_2$ of exact period $n>1$.
\begin{Lemma}
\label{Lem:LutzNagell}
The following statements hold.
\begin{itemize}
\item[a)] If $p$ is a prime such that $\nu_p( u_1)=a$, then $\nu_p(F_0)=a$. In particular, $u_1\parallel F_0$.
\item[b)] $c_1$ and $u_1$ are relatively prime in $O_K$.
\item[c)] $\displaystyle u_1\parallel \frac{F_0}{c_1}$; and $\displaystyle\frac{F_0}{c_1}$ and $c_1$ are relatively prime in $O_K$.
\item[d)] $\displaystyle c_1\mid (u_1^{d^{n}-1}-u_2^{d^{n}-1})$.
\end{itemize}
\end{Lemma}
\begin{Proof}
We will be mainly considering eq (\ref{eq1}) above. For (a), that $\nu_p( F_0)\ge a$ is a direct consequence of eq (\ref{eq1}) and the fact that $u_1$ and $u_2$ are relatively prime. If $p^{a+1}\mid F_0$, then this will imply that $p$ divides the coefficient of the linear term in $u_1$, namely, $c_2^{d^{n-1}}u_2^{d^n-1}$, which is a contradiction.

For (b), according to Corollary \ref{cor:conjecture}, the linear factor $\displaystyle \left(x-(g\circ h)^j(u_1)/u_2\right)$, $1\le j\le n$, divides $G^n(x)$. In other words, $u_2x-(g\circ h)^j(u_1)$ divides $u_2^{d^n} G^n(x/u_2)$. In particular, one sees that $u_1(u_1^d+c_1)/u_2^{d-1}$ divides $F_0$. It follows that if there is a common prime divisor $p$ of $c_1$ and $u_1$ such that $\nu_p( u_1)=a$, then $\nu( F_0)>a$ which contradicts (a).

Since $\displaystyle F_0= c_1c_2^{d^{n-1}-1}+\sum_{i=0}^{d^{n-1}-d}g_{n,i}c_1^{d+i}c_2^{d^{n-1}-d-i}\in c_1\Z[c_1,c_2]$, see Lemma \ref{lem:freecoefficient}, part (c) follows directly from (a) and (b) and the condition that $c_1$ and $c_2$ are relatively prime in $O_K$.

Since $F_i\in c_1\Z[c_1,c_2]$, $i\ne d^{n-1}$, it follows that \[c_1\mid F_{d^{n-1}}u_1^{d^n}-c_2^{d^{n-1}}u_1u_2^{d^n-1}= c_2^{d^{n-1}}u_1(u_1^{d^{n}-1}-u_2^{d^{n}-1}).\] Since $c_1$ is relatively prime to both $u_1$ and $u_2$ in $O_K$, where the latter relative primality holds because $c_2=u_2^d$, this yields that $c_1\mid (u_1^{d^{n}-1}-u_2^{d^{n}-1})$.
\end{Proof}

\section{Periodic points and divisors of arithmetic sequences}
In the rest of this note, we illustrate the connection between periodic points of the polynomial $f(x)=x^d+c\in \Q[x]$ and two arithmetic sequences.

Let $c=c_1/c_2$ be such that $c_1\in\Z$ and $c_2\in\Z^+$ are relatively prime. Given that $u_1/u_2$ is a periodic point of exact period $n$ of $x^d+c$, the orbit of $u_1/u_2$ is the set $\displaystyle O_f(u_1/u_2)=\{f^{j}(u_1/u_2):j=1,2,3,\ldots\}$. We recall that $f^j(u_1/u_2)=(g\circ h)^j(u_1)/u_2$ where $h(x)=x^d+c_1$ and $g(x)=x/u_2^{d-1}$, $j=1,2,\ldots$, see Remark \ref{rem1} and Corollary \ref{cor:conjecture}. We set $u_{1,j}=(g\circ h)^j(u_1)$.

In this section, fixing $i$ and $j$, we consider the sequence $\displaystyle \frac{u_{1,i}^k-u_{1,j}^{k}}{u_{1,i}-u_{1,j}}$, $k=1,2,3,\ldots$. We investigate the divisibility of the terms of the latter sequence by prime divisors of $c_1$. In fact, according to Lemma \ref{Lem:LutzNagell} d), if $p$ is a prime divisor of $c_1$, then $p\mid  \left(u_{1,l}^{d^n-1}-u_2^{d^n-1}\right)$ for every $l$. Therefore, $p\mid  (u_{1,i}^{d^n-1}-u_{1,j}^{d^n-1})$ for any $i$ and $j$.

We first prove the coprimality of $u_{1,i}$ and $u_{1,j}$ for any choice of $i$ and $j$, $i\ne j$.

\begin{Lemma}
\label{lem:coprimality}
Let $f(x)=x^d+c_1/c_2\in K[x]$ where $c_1\in O_K$ and $c_2\in O_K/O_K^{\times}$ are relatively prime. If $u_1/u_2$ is a periodic point of exact period $n$, where $u_1$ and $u_2$ are relatively prime in $O_K$, then $u_{1,i}$ and $u_{1,j}$ are relatively prime for any $i\ne j$.
\end{Lemma}
\begin{Proof}
Let $p$ be a common prime divisor of $u_{1,i}$ and $u_{1,j}$. Assume that $\nu_p (u_{1,k})=a_k$, $k=i,j$. According to Lemma \ref{Lem:LutzNagell}, one has $\nu_p(F_0)=a_i=a_j$ where $F_0$ is defined as before. Since both $u_{1,i}/u_2$ and $u_{1,j}/u_2$ are periodic points of $f(x)$, it follows that they are zeros of the polynomial $G^n(x)$ defined in \S 4. In particular, $u_{1,i}u_{1,j}$ divides $F_0$. Therefore, if $p$ was a prime divisor of both $u_{1,i}$ and $u_{1,j}$, this would contradict the fact that $\nu_p(F_0)=a_i$.
\end{Proof}
\begin{Theorem}
\label{thm:primitive}
  Let $u_1/u_2$ be a periodic point of $f(x)=x^d+c\in\Q[x]$ of exact period $n$ where $c=c_1/c_2$ is as above. Assume, moreover, that there is a prime $p\mid c_1$ such that $\gcd(p,d^n-1)=1$, then $\displaystyle p\nmid (u_{1,i}-u_{1,j})$, for all $i\ne j$.
  In particular, $\displaystyle p\mid \frac{u_{1,i}^{d^n-1}-u_{1,j}^{d^n-1}}{u_{1,i}-u_{1,j}}$.
\end{Theorem}
\begin{Proof}
Let $p$ be a prime such that $p|c_1$ and $\gcd(p,d^n-1)=1$. We assume on the contrary that $\nu_p (u_{1,i}-u_{1,j})=\alpha>0$. We set $\displaystyle b_{i,j}(m)=\frac{u_{1,i}^m-u_{1,j}^m}{u_{1,i}-u_{1,j}}$. We recall that
\begin{eqnarray*}
\gcd(b_{i,j}(k),b_{i,j}(l))=b_{i,j}(g),\qquad g=\gcd(k,l),
\end{eqnarray*}
see \cite[Theorem VI]{Carmichael}.

 Since $\nu_p(u_{1,i}-u_{1,j})=\alpha$, one has $\nu_p(u_{1,i}^p-u_{1,j}^p)\ge \alpha+1$, see \cite[Theorem III]{Birkhoff}. Noting that $\displaystyle \gcd\left(b_{i,j}(m),b_{i,j}(p)\right)=b_{i,j}(1)=1$ whenever $\gcd(m,p)=1$ and that $\nu_p(u_{1,i}^k-u_{1,j}^k)\ge \alpha$ for all $k\ge 1$, one has $\nu_p(u_{1,i}^m-u_{1,j}^m)=\nu_p(u_{1,i}-u_{1,j})=\alpha$ whenever $\gcd(m,p)=1$.

Since $u_{1,i}/u_2$ is a point in the orbit of $u_1/u_2$, hence a periodic point of period $n$, one has $f^{n}(u_{1,i}/u_2)=u_{1,i}/u_2$. Thus, eq (\ref{eq1}) may be written for $u_{1,i}/u_2$ as follows
\begin{eqnarray}
\label{eq11}
F_0u_2^{d^{n}}+F_1u_{1,i}^du_2^{d^{n}-d}+F_2u_{1,i}^{2d}u_2^{d^{n}-2d}+\ldots+F_{d^{n-1}-1}u_{1,i}^{d^{n}-d}u_2^d+F_{d^{n-1}}u_{1,i}^{d^{n}}=c_2^{d^{n-1}}u_{1,i}u_2^{d^{n}-1}.\nonumber\\
\end{eqnarray}
Similarly,
\begin{eqnarray}
\label{eq12}
F_0u_2^{d^{n}}+F_1u_{1,j}^du_2^{d^{n}-d}+F_2u_{1,j}^{2d}u_2^{d^{n}-2d}+\ldots+F_{d^{n-1}-1}u_{1,j}^{d^{n}-d}u_2^d+F_{d^{n-1}}u_{1,j}^{d^{n}}
=c_2^{d^{n-1}}u_{1,j}u_2^{d^{n}-1}.\nonumber\\
\end{eqnarray}
Multiplying (\ref{eq11}) and (\ref{eq12}) times $u_{1,j}^{d^{n}}$ and $u_{1,i}^{d^{n}}$, respectively, and subtracting the two resulting equations, one obtains
{\footnotesize\begin{align}\label{eq3}F_0u_2^{d^n}(u_{1,i}^{d^n}-u_{1,j}^{d^n})+F_1\left(u_{1,i}^{d^n-d}-u_{1,j}^{d^n-d}\right)u_{1,i}^du_{1,j}^du_2^{d^{n}-d}+F_2\left(u_{1,i}^{d^n-2d}-u_{1,j}^{d^n-2d}\right)u_{1,i}^{2d}u_{1,j}^{2d}u_2^{d^{n}-2d}
+\ldots\nonumber\\+F_{d^{n-1}-1}\left(u_{1,i}^{d}-u_{1,j}^{d}\right)u_{1,i}^{d^n-d}u_{1,j}^{d^n-d}u_2^d=c_2^{d^{n-1}}\left(u_{1,i}^{d^n-1}-u_{1,j}^{d^n-1}\right)u_{1,i}u_{1,j}u_2^{d^{n}-1}.
\end{align}}
 One recalls that $F_i\in c_1\Z[c_1,c_2]$ for $i\ne d^{n-1}$, see Corollary \ref{cor1}, and $p^{\alpha}||(u_{1,i}-u_{1,j})$. This yields that the left hand side of eq (\ref{eq3}) is divisible by $p^{\alpha+1}$. Now since $c_1$ is relatively prime to each of $c_2$, $u_2$, $u_{1,i}$ and $u_{1,j}$, it follows that $p^{\alpha+1}$ divides $\displaystyle\left(u_{1,i}^{d^n-1}-u_{1,j}^{d^n-1}\right)$ on the right hand side of eq (\ref{eq3}), which is a contradiction as $\gcd(p,d^n-1)=1$.
\end{Proof}

\begin{Corollary}
\label{cor:primitive1}
Let $u_1/u_2$ be a periodic point of $x^d+c$ of exact period $n$ where $c=c_1/c_2$ is as above. If there is a prime $p$ such that $p\mid c_1$ and $\gcd(p, d^n-1)=1$, then $\gcd(p-1,d^n-1)>1$. In fact, if $d^n-1$ is prime, then $p\equiv 1$ mod $(d^n-1)$, in particular, $p>d^n$.
\end{Corollary}
\begin{Proof}
Since $\gcd(p,d^n-1)=1$, one knows that $p\nmid (u_{1,i}-u_{1,j})$, see Theorem \ref{thm:primitive}. We recall that
\begin{eqnarray*}
\gcd(b_{i,j}(k),b_{i,j}(l))=b_{i,j}(g),\qquad g=\gcd(k,l).
\end{eqnarray*}
Since $\nu_p\left(u_{1,i}^{p-1}-u_{1,j}^{p-1}\right)>0$ by Fermat's Little Theorem, one knows that $\nu_p(b_{i,j}(p-1))>0$. Furthermore, as $c_1\mid\left(u_{1,i}^{d^n-1}-u_{1,j}^{d^n-1}\right)$, one has $\nu_p(b_{i,j}(d^n-1))>0$. It follows that $\gcd(p-1, d^n-1)>1$.

If $d^{n}-1$ is prime, then $d^n-1$ is the order of $u_1u_2^{-1}$ mod $p$. This implies that $(d^n-1)\mid p-1$.
\end{Proof}

\begin{Remark}
Let $p$ be a prime divisor of $c_1$ such that $\gcd(p,d^n-1)=1$. In view of Corollary \ref{cor:primitive1}, if $\gcd(p-1,d^{n}-1)=1$, then $x^d+c_1/c_2$ has no periodic points of period $n$. Furthermore, if $d^n-1$ is prime, then $d^n-1$ divides $p-1$ for every prime divisor $p$ of $c_1$.
Finally, if $p\mid (u_{1,i}^m-u_{1,j}^m)$ for some $m<(d^n-1)$, then $\gcd(m,p-1)>1$. In particular, if $\gcd(m,p-1)=1$ for any $m<d^n-1$, then $p$ is a primitive prime divisor of $\displaystyle \frac{u_{1,i}^{d^n-1}-u_{1,j}^{d^n-1}}{u_{1,i}-u_{1,j}}$.
\end{Remark}
\begin{Example}
Let $m> 1$. Let the polynomial $f(x)=x^2+2^m$ be such that $2^m-1$ is prime. If $n>1$ is an integer such that $\gcd(m,n)=1$, then $\gcd(2^m-1,2^n-1)=1$. Thus, Corollary \ref{cor:primitive1} implies that $f(x)=x^2+2^m$ has no periodic point of period $n$ when $\gcd(m,n)=1$.
\end{Example}

\section{A remark on primitive prime divisors of $f^n(0)$}
We recall that if $x_i,i=1,2,\ldots$, is a sequence in the ring of integers $O_K$ of a number field $K$, then the term $x_n$ is said to have a {\em primitive prime divisor} $p$ if $p$ is a prime such that $\nu_p(x_n)>0$, and $\nu_p(x_m)=0$ for any $m<n$.

Set $f(x)=x^d+c_1/c_2\in K[x]$, $c_1\in O_K$, $c_2\in O_K/O_K^{\times}$, $d\ge 2$. In this section, we write $F_0^n$ for $c_2^{d^{n-1}}f^n(0)$. It is known that the sequence $F_0^n$ is a divisibility sequence. In particular, $F_0^m\mid F_0^n$ whenever $m\mid n$. Several results were proved concerning the existence of primitive prime divisors for each term of the sequence $F_0^n$, see for example \cite{Krieger}.
\begin{Lemma}
\label{lemma1}
Let $K$ be a number field with ring of integers $O_K$. Let $g(x)\in O_K[x]$ and $u\in O_K$ be such that there is a prime $p$ dividing $g^m(u)$ and $g^n(u)$, $ n>m$. Then $p$ divides $g^{n-m}(0)$.
\end{Lemma}
\begin{Proof}
This follows directly by observing that $g^n(u)=g^{n-m}(g^m(u))$.
\end{Proof}

\begin{Theorem}
\label{thm:primitiveF0}
If $u_1/u_2$ is a periodic point of $f(x)=x^d+c_1/c_2\in K[x]$ of exact period $n$, where $u_i,c_i$ are as before, then every prime divisor of $u_1$ is a primitive prime divisor of $F_0^n$, $n>1$.
\end{Theorem}
\begin{Proof}
 One knows that $u_1\mid (F_0^n/c_1)$, see Lemma \ref{Lem:LutzNagell} c). We assume that $p$ is a prime divisor of $u_1$ such that $p\mid F_0^m$ for $m<n$. According to Lemma \ref{lemma1}, one has $\nu_p( F_0^{n-m})>0$. Let $m$ be the smallest such positive integer. One knows that $m\ge2$ since $\gcd(c_1,u_1)=1$, see Lemma \ref{Lem:LutzNagell} b). By successive application of the division algorithm, one has $m\mid n$.

  Therefore, if $n$ is prime, then it is impossible for $p$ to divide $F_0^m$ for $m<n$.

   Now, we assume $n$ is composite. Let $q_1$ and $q_2$ be two distinct prime divisors of $n$ where $n=q_ik_i$. We consider the polynomial $g_i(x)=f^{k_i}(x)$. One has $g_i(0), g_i^2(0)=f^{2k_i}(0)$, $g_i^3(0)=f^{3k_i}(0),\ldots,g_i^{q_i}(0)=f^n(0)$. Since $f^{n}(0)=g_i^{q_i}(0)$, Lemma \ref{Lem:LutzNagell} implies that $\nu_p(g_i^{q_i}(0))>0$. Since $q_i$ is prime, it follows that the smaller possible integer $l$ such that $\nu_p(g_i^l(0))>0$ is $l=1$. In other words, $\nu_p(f^{k_1}(0)),\nu_p(f^{k_2}(0))>0$. This yields that either $k_1\mid k_2$ or $k_2\mid k_1$, a contradiction.
 \end{Proof}

\begin{Corollary}
If $f(x)=x^d+c_1/c_2\in\Q[x]$ has a periodic point of period $n$, then $F_0^n$ has at least $n-1$ distinct primitive prime divisors.
\end{Corollary}
\begin{Proof}
This follows immediately from Theorem \ref{thm:primitiveF0} and Lemma \ref{lem:coprimality}.
\end{Proof}

\end{document}